\documentclass[11pt,dvips,twoside,letterpaper]{article}
\usepackage{pslatex}
\usepackage{fancyhdr}
\usepackage{graphicx}
\usepackage{geometry}

\usepackage{amsmath}
\usepackage{amssymb}
\usepackage{amsfonts}
\usepackage{amsthm,amscd}

\newtheorem{theorem}{Theorem}[section]

\newtheorem{corollary}[theorem]{Corollary}
\newtheorem{proposition}[theorem]{Proposition}
\theoremstyle{definition}
\newtheorem{definition}[theorem]{Definition}
\newtheorem{example}[theorem]{Example}

\theoremstyle{remark}
\newtheorem{remark}[theorem]{Remark}

\numberwithin{equation}{section}

\def\C{\mathbb C}

\def\R{\mathbb R}

\def\Z{\mathbb Z}

\def\N{\mathbb N}

\newcommand{\eps}{\varepsilon}
\newcommand{\al}{\mathbb{\alpha}}
\newcommand{\diam}{\diamondsuit_{\alpha}}

\newcommand{\rar}{\mbox{$\rightarrow$}}

\newcommand{\T}{\mathbb{T}}

\newcommand{\ps}{p^{\sigma}}
\newcommand{\expal}{\mbox{}_{\alpha}e_p}

\begin{document}

\title{\bfseries\scshape{A Study of Diamond-alpha dynamic equations on regular time scales}\thanks{Accepted (09/February/2009) for publication
at the \emph{African Diaspora Journal of Mathematics},
http://www.african-j-math.org}}

\author{\bfseries\scshape Dorota Mozyrska\thanks{E-mail address: admoz@w.tkb.pl}\\
Faculty of Computer Science, Bia{\l}ystok Technical
University\\15-351 Bia\l ystok, Poland\\
\\\bfseries\scshape Delfim F. M. Torres\thanks{On leave of absence from Department of Mathematics, University of Aveiro,
3810-193 Aveiro, Portugal. E-mail address: delfim@ua.pt}\\
African University of Science and Technology\\
Galadimawa P.M.B 681, Garki Abuja, Nigeria}

\date{}
\maketitle
\thispagestyle{empty}
\setcounter{page}{1}


\thispagestyle{fancy} \fancyhead{}
\fancyhead[L]{}
\fancyfoot{}
\renewcommand{\headrulewidth}{0pt}


\begin{abstract}
We introduce the diamond-alpha exponential function on time scales. As particular cases, one gets both delta and nabla exponential functions. A method of solution of a homogenous linear dynamic diamond-alpha equation on a regular time scale is investigated, and examples of diamond-alpha exponential functions are presented.
\end{abstract}


\vspace{.08in} \noindent \textbf{Keywords}: Time scales, exponential function, diamond-alpha derivatives.


\section{Introduction}

The calculus on time scales is a powerful tool to unify discrete and continuous analysis \cite{AB,Bh,Bh1}. The generalization to time scales of forward and backward differences in $\mathbb{Z}$ leads to delta and nabla derivatives, respectively. Based on these two types of derivatives, a combined derivative, the so called diamond-$\alpha$ derivative, was recently introduced by Sheng, Fadag, Henderson, and Davis \cite{SFHD}. For $\alpha=1$ the diamond-$\alpha$ derivative reduces to the delta-derivative while for $\alpha=0$ we get the nabla-derivative. When one chooses the time scale to be the set of real numbers, then all such notions are equivalent and coincident with the standard concept of derivative. For other time scales, the diamond-$\alpha$ derivative is a linear convex combination of delta and nabla-derivatives and seems to be useful for solving certain dynamic equations \cite{Sheng1,SFHD}. Although the combined diamond-$\alpha$ derivative has no anti-derivative, and so it is not a dynamic derivative \cite{Sheng1}, it is still possible to define a diamond-$\alpha$ integral on time scales with many interesting properties \cite{diamondBasia,jia,gruss}.
The diamond-$\alpha$ calculus
is still in its infancy, and much remains
to be done. Here we introduce the diamond-$\alpha$ exponential function from an appropriate diamond-$\alpha$ eigenfunction.

It is known that eigenvalue problems of delta-differential operators play an important role on time scales \cite{SLP,iso:ts}. We look to the diamond-$\alpha$ derivative as a linear operator on the linear space of differentiable functions in both delta and nabla senses. More precisely, we investigate and compare the idea of eigenfunction on time scales for three kinds of differential operators on a space of functions. An eigenfunction of a linear operator $L$, defined on some function space, is understood as any non-zero function $f$ in that space for which there exists a function $\lambda$ such that $Lf=\lambda f$. The solution of such problems depends on boundary conditions required for $f$.  Typically, the exponential functions are defined on time scales by means of a cylindrical mapping \cite{Bh,Bh1}. An alternative but equivalent way, is to define them as the solution of an eigenvalue problem. A suitable $(Ly)(\cdot)$  gives the exponential function as the solution of a problem $Ly=py$, where $p$ is a function, subject to appropriate boundary conditions.


\pagestyle{fancy}
\fancyhead{}
\fancyhead[EC]{D. Mozyrska and D. F. M. Torres}
\fancyhead[EL,OR]{\thepage} \fancyhead[OC]{A Study
     of Diamond-alpha Dynamic Equations on Regular Time Scales}
\fancyfoot{}
\renewcommand\headrulewidth{0.5pt}


\section{Preliminaries}

In this section we give a short introduction, with basic definitions and necessary results to what follows, on three types of calculus on time scales: (i) the delta ($\Delta$),
(ii) the nabla ($\nabla$), and (iii) the diamond-$\alpha$ ($\diam$) calculus. For more information we refer the reader to \cite{AB,Bh}, \cite{And,Bh}, and \cite{Ozkan,Sheng1,Sheng,SFHD}, respectively.

By a time scale, denoted by $\T$, we mean a nonempty closed
subset of $\R$. As the theory of time scales give a way to unify
continuous and discrete analysis, the standard cases of time
scales are $\T=\R$, $\T=\Z$, $\T=\N$ or $\T=c\Z$, $c>0$.
In particular also the set $Q_q=\{-q^k,0,q^k, k\in\Z, q>1\}$ forms a time scale. The set of natural numbers is taken here with zero: $\N=\{0, 1, 2, \ldots\}$.
Let $\T^*:=\T\backslash\{\min\T, \max \T\}$. Then $\N^*=\{1, 2, \ldots\}$, while for example, $c\Z^*=c\Z$.

For $t\in\T$, the forward jump operator $\sigma$ and the
graininess function $\mu$ are defined by
$\sigma(t)=\inf\{s\in\T:s>t\}$ and $\sigma(\sup\T)=\sup\T$ if
$\sup\T<+\infty$; $\mu(t)=\sigma(t)-t$.  Moreover,
we define the backward
operator $\rho$ by $\rho(t)=\sup\{s\in\T: s<t\}$ and
$\rho(\inf\T)=\inf\T$ if $\inf\T>-\infty$; $\nu(t)= t - \rho(t)$.
In the continuous-time
case, \textrm{i.e.}, when $\T=\R$, we have $\sigma (t)=\rho(t)=t$
and $\mu(t)=\nu(t)=0$ for all $t\in\R$.
In the discrete-time case,
$\sigma(t)=t+1$, $\rho(t)=t-1$, and $\mu(t)=\nu(t)=1$ for each
$t\in\T=\Z$. For the composition between a function
$f:\T\rightarrow\R$ and functions $\sigma:\T\rightarrow\T$ and
$\rho:\T\rightarrow\T$, we use the abbreviations
$f^{\sigma}(t)=f(\sigma(t))$ and $f^{\rho}(t)=f(\rho(t))$.
A point $t$ is called
left-scattered (right-scattered) if $\rho(t)<t$ ($\sigma(t)>t$). A point $t$ is called left-dense (right-dense) if
$\rho(t)=t$ ($\sigma(t)=t$).

The set $\T^{\kappa}$ is defined by $\T^{\kappa}:=\T\backslash (\rho(\sup\T),\sup\T]$ if
$\sup\T<\infty$, and $\T^{\kappa}=\T$ if $\sup\T=\infty$; the set $\T_{\kappa}$
by $\T_{\kappa}:=\T\backslash[\inf\T,\sigma(\inf\T))$ if $|\inf\T|
<\infty$, and $\T_{\kappa}=\T$ if $\inf\T=-\infty$. Moreover, $\T^{\kappa^{n+1}}:=\left(\T^{\kappa^n}\right)^\kappa$,
$\T_{\kappa^{n+1}}:=\left(\T_{\kappa^n}\right)_{\kappa}$ and
$\T^{\kappa}_{\kappa}:=\T^{\kappa}\cap \T_{\kappa}$.

For a function $f:\T\rightarrow\R$, we define the
$\Delta$--derivative of $f$ at $t\in\T^{\kappa}$, denoted by
$f^{\Delta}(t)$, to be the number, if it exists, with the property
that for all $\varepsilon>0$, exists a neighborhood $U\subset \T$
of $t\in\T^{\kappa}$ such that for all $s\in U$,
$|f^\sigma(t)-f(s)-f^{\Delta}(t)(\sigma(t)-s)|\leq
\varepsilon|\sigma(t)-s|$. Function $f$ is said to be
$\Delta$--differentiable on $\T^{\kappa}$
provided $f^{\Delta}(t)$ exists
for all $t\in\T^{\kappa}$.

The $\nabla$--derivative of $f$, denoted by $f^{\nabla}(t)$, is
defined in a similar way: it is the number, if it exists, such
that for all $\varepsilon>0$ there is a neighborhood $V\subset \T$
of $t\in\T_{\kappa}$ such that for all $s\in V$,
$|f^\rho(t)-f(s)-f^{\nabla}(t)(\rho(t)-s)|\leq
\varepsilon|\rho(t)-s|$. Function $f$ is said to be
$\nabla$--differentiable on $\T_{\kappa}$ provided $f^{\nabla}(t)$ exists for all $t\in\T_{\kappa}$.

\begin{example} The classical settings are obtained
by choosing $\T = \R$ or $\T=c\Z$, $c>0$:
\begin{enumerate}
\item Let $\T=\R$ and $f$ be differentiable.
Then, $f^{\Delta}(t)=f^{\nabla}(t)=f'(t)$ and $f$ is $\Delta$ and $\nabla$ differentiable if and only if it is differentiable in the ordinary sense.

\item Let $\T=c\Z$. Then, both derivatives
$$
f^{\Delta}(t)=\frac{f(t+c)-f(t)}{c} \, , \quad
f^{\nabla}(t)=\frac{1}{c}\left(f(t)-f(t-c)\right)
$$
always exist, $t \in \T$.
\end{enumerate}
\end{example}

It is possible to establish some relationships between $\Delta$
and $\nabla$ derivatives.

\begin{theorem}{\rm \cite{Bh}}
\label{th:1}
(a) Assume that $f:\T\rightarrow\R$ is $\Delta$--differentiable
on $\T^{\kappa}$. Then, $f$ is $\nabla$--differentiable at $t$ and $f^{\nabla}(t)=f^{\Delta}(\rho(t))$ for all $t\in\T_{\kappa}$ such that $\sigma(\rho(t))=t$. (b) Assume that $f:\T\rightarrow\R$ is $\nabla$--differentiable on $\T_{\kappa}$. Then, $f$ is $\Delta$--differentiable at $t$
and $f^{\Delta}(t)=f^{\nabla}(\sigma(t))$ for all $t\in\T^{\kappa}$ such
that $\rho(\sigma(t))=t$.
\end{theorem}

A function $f:\T \rar \R$ is called rd-continuous provided it is
continuous at right-dense points in $\T$ and its left-sided limits exist (finite) at left-dense points in $\T$. The class of real
rd-continuous functions defined on a time scale $\T$ is denoted by $C_{rd}(\T,\R)$. If $f\in C_{rd}(\T,\R)$, then there exists a function $F(t)$ such that $F^{\Delta}(t)=f(t)$. The $\Delta$--integral is defined by $\int_{a}^bf(t)\Delta t=F(b)-F(a)$.

Similarly, a function $f:\T \rar \R$ is called ld-continuous
provided it is continuous at left-dense points in $\T$ and its
right-sided limits exist (finite) at right-dense points in $\T$.
The class of real ld-continuous functions defined on a time scale
$\T$ is denoted by $C_{ld}(\T,\R)$. If $f\in C_{ld}(\T,\R)$, then
there exists a function $G(t)$ such that $G^{\nabla}(t)=f(t)$. In
this case we define $\int_{a}^bf(t)\nabla t=G(b)-G(a)$.

\begin{definition}{\rm \cite{Sheng1}}
Let $\T$ be a time scale, $\al\in[0,1]$,
$\mu_{ts}=\sigma(t)-s$, $\eta_{ts}=\rho(t)-s$, and
$f:\T\rightarrow\R$. The $\diam$--derivative of $f$ at $t$ is
defined to be the value $f^{\diam}(t)$, if it exists, such that
for all $\varepsilon>0$ there is a neighborhood $U\subset \T$ of
$t$ such that for all $s\in U$,
$$|\al\left[f^{\sigma}(t)-f(s)\right]\eta_{ts}
+(1-\al)\left[f^{\rho}(t)-f(s)\right]\mu_{ts}
-f^{\diam}(t)\mu_{ts}\eta_{ts}|\leq
\varepsilon|\mu_{ts}\eta_{ts}| \, .$$
We say that function $f$ is
$\diam$--differentiable on $\T^{\kappa}_{\kappa}$,
provided $f^{\diam}(t)$ exists
for all $t\in\T^{\kappa}_{\kappa}$.
\end{definition}

\begin{theorem}{\rm \cite{Sheng1}}
\label{thsheng}
Let $f:\T\rightarrow\R$  be simultaneously $\Delta$ and $\nabla$
differentiable at $t\in\T$. Then, $f$ is $\diam$--differentiable at $t$ and $f^{\diam}(t)=\al f^{\Delta}(t)+(1-\al)f^{\nabla}(t)$,
$\al\in[0,1]$.
\end{theorem}

\begin{remark} If a given function is $\Delta$ and $\nabla$ differentiable at the point $t$, then the $\diam$--derivative at $t$ is the convex combination of $\Delta$ and $\nabla$
derivatives. It reduces to the $\Delta$--derivative for $\al=1$ and to the $\nabla$--derivative for $\al=0$. The case $\al=0.5$ has proved to be very useful in applications \cite{SFHD,Sheng1,Sheng,Ozkan}.
\end{remark}

We use the following notation:
$C_{rl}(\T,\R):=C_{rd}(\T,\R)\cap C_{ld}(\T,\R)$.

\begin{definition}\label{diaminteg}
Let $a,b\in\T$ and $f\in C_{rl}(\T,\R)$. Then,
the $\diam$--integral of $f$ is defined by
$\int_a^bf(\tau)\diam\tau=\al
\int_a^bf(\tau)\Delta\tau+(1-\al)\int_a^bf(\tau)\nabla\tau$, where
$\al\in[0,1]$.
\end{definition}

It should be noted that, in general, the $\diam$--derivative of $\int_a^t f(\tau)\diam\tau$ with respect to $t$ is not equal to $f(t)$.


\section{Delta and Nabla Exponential Functions}

A function $p:\T\rightarrow \R$ is called regressive provided $1+\mu(t)p(t)\neq 0$ for all $t\in\T^{\kappa}$.
By $ \mathcal{R}$ it is denoted the set of all regressive and rd-continuous functions on $\T$.
Similarly, a function $q:\T\rightarrow \R$ is called $\nu$--regressive provided $1-\nu(t)q(t)\neq 0$ for all $t\in\T_{\kappa}$.
By $ \mathcal{R}_{\nu}$ we denote the set of all $\nu$--regressive and ld-continuous functions on $\T$.
We recall now some notations about Hilger's complex plane \cite{AP,Bh}. By $\C_h$ it is denoted the set of Hilger's complex numbers: $\C_0:=\C$, $\C_h:=\{z\in\C:z\neq-\frac{1}{h}\}$ for $h>0$; whereas by $\R_h$ we denote Hilger's real axis:  $\R_0:=\R$, $\R_h:=\{z\in\C_h: z\in\R \ \mbox{and}\ z>-\frac{1}{h}\}$ for $h > 0$.
Let $\Z_h$ be the strip $\Z_h:=\{z\in\C:-\frac{\pi}{h}< Im (z)<\frac{\pi}{h}\}$, $h>0$. Then, one defines the cylinder and $\nu$--cylinder transformations as follows:

a) $\xi_h:\C_h\rightarrow\Z_h$ by $\xi_h(z):=\frac{1}{h}Log(1+zh)$,

b) $\hat{\xi}_h:\C_h\rightarrow\Z_h$ by $\hat{\xi}_h(z):=-\frac{1}{h}Log(1-zh)$,

\noindent where $Log$ is the principal logarithm function \cite{AP,Bh}. For $h=0$, one has $\xi_0(z)=z$ and $\hat{\xi}_0(z)=z$ for all $ z\in\C$.

\begin{definition}\mbox{}\\
i) Let
$p\in \mathcal{R}$. Then, the $\Delta$--exponential function is defined by
\begin{equation}\label{exp1}
e_p(t,t_0):=\exp\left(\int_{t_0}^t\xi_{\mu(\tau)}(p(\tau))\Delta\tau\right).
\end{equation}
ii) Let $p\in \mathcal{R}_{\nu}$. Then, the $\nabla$--exponential function is defined by
\begin{equation}\label{exp2}
\hat{e}_p(t,t_0):=\exp\left(\int_{t_0}^t\hat{\xi}_{\nu(\tau)}(p(\tau))\nabla\tau\right).
\end{equation}
\end{definition}

\begin{proposition}\mbox{\cite{AP,Bh}}\\
a) Let
$p\in \mathcal{R}$. Then, the $\Delta$--exponential function $e_p(\cdot,t_0)$ is the unique solution of the initial value problem
$y^{\Delta}(t)=p(t)y(t)$, $y(t_0)=1$.\\
b) Let $p\in \mathcal{R}_{\nu}$. Then, the $\nabla$--exponential function $\hat{e}_p(\cdot,t_0)$ is the unique solution of the initial value problem
$y^{\nabla}(t)=p(t)y(t)$, $y(t_0)=1$.
\end{proposition}

The next theorem establishes a relation between the
$\Delta$ and $\nabla$ exponential functions.

\begin{theorem}{\rm \cite[Theorem~4.17]{Bh1}}
If $p$ is continuous and regressive, then
\begin{equation*} e_p(t,t_0)=\hat{e}_{\frac{p^{\rho}}{1+p^{\rho}\nu}}(t,t_0).
\end{equation*}
If $q$ is continuous and $\nu$--regressive, then
\begin{equation*} \hat{e}_q(t,t_0)=e_{\frac{q^{\sigma}}{1-q^{\sigma}\mu}}(t,t_0).
\end{equation*}
\end{theorem}

We need also the following:

\begin{corollary}\label{cor:20.01}
Let $p(\cdot)\in \mathcal{R}\cap \mathcal{R}_{\nu}$. Then,

a) $e_p^{\rho}(t,t_0)=\frac{1}{1+p^{\rho}(t)\nu(t)}e_p(t,t_0)$, and
$e_p^{\rho}(t_0,t_0)=\frac{1}{1+p^{\rho}(t_0)\nu(t_0)}$;

b) $\hat{e}_p^{\rho}(t,t_0)=1-p(t)\nu(t)\hat{e}_p(t,t_0)$,
and
$\hat{e}_p^{\rho}(t_0,t_0)=1-p(t_0)\nu(t_0)$.
\end{corollary}


\section{Naive Approach: Combined-Exponentials}

Roughly speaking, the $\diam$--calculus is the convex combination of $\Delta$ and $\nabla$ calculuses. We may be then tempted
to define the $\diam$--exponential function by a simple combination of $\Delta$ and $\nabla$ exponentials. We consider here two such functions: $\mbox{}_{\alpha}E_p$ and $\mbox{}_{\alpha}e_p$,
where $p\in \mathcal{R}\cap \mathcal{R}_{\nu}$ and $\al\in[0,1]$.

Our first combined-exponential is the most natural to be considered:
\begin{equation*}
\mbox{}_{\alpha}E_p(\cdot,t_0) := \al e_p(\cdot, t_0)+(1-\al)\hat{e}_p(\cdot,t_0), \quad t_0\in\T.
\end{equation*}

\begin{example}
Consider the time scale $\T=\Z$ and the constant function $p(t)=\frac{1}{2}$. Take $t_0=0$. Then, $e_p(t,0)=\left(\frac{3}{2}\right)^t$ is the solution of the initial value problem $y^{\Delta}(t)=\frac{1}{2}y(t)$, $y(0)=1$.  Moreover, $\hat{e}_p(t,0)=2^t$ is the unique solution of $y^{\nabla}(t)=\frac{1}{2}y(t)$, $y(0)=1$. It follows that
$\mbox{}_{\alpha}E_p(t,0)
=\al\left(\frac{3}{2}\right)^t+(1-\al)2^t$, $t\in\Z$.
\end{example}

Our second combined-exponential $\mbox{}_{\alpha}e_p$
also involves a convex combination, but of a different nature.
It is motivated by the definitions \eqref{exp1} and \eqref{exp2}
of $\Delta$-- and $\nabla$--exponential functions:
\begin{equation*}
\mbox{}_{\alpha}e_p(t,t_0)
:=\exp\left(\alpha\int_{t_0}^t\xi_{\mu(\tau)}(p(\tau))\Delta\tau+
(1-\alpha)\int_{t_0}^t\hat{\xi}_{\nu(\tau)}(p(\tau))\nabla\tau\right),
\quad  t_0, t \in \T.
\end{equation*}
The combined-exponential $\mbox{}_{\alpha}e_p$ presents some desired properties:

\noindent (i)  $\mbox{}_{\alpha}e_p(t,t_0)=e_p^{\al}(t,t_0)\hat{e}_p^{1-\al}(t,t_0)$;

\noindent (ii) $\ln\left(\mbox{}_{\alpha}e_p(t,t_0)\right)=\al e_p(t,t_0)+(1-\alpha)\hat{e}_p(t,t_0)$;

\noindent (iii) if $p\in \mathcal{R}\cap \mathcal{R}_{\nu}$, then the semigroup property $\expal (t,s)\expal(s,t_0)=\expal(t,t_0)$
holds.

While both functions $\mbox{}_{\alpha}E_p$ and $\mbox{}_{\alpha}e_p$
generalize $\Delta$ and $\nabla$ exponentials, in the sense that when $\alpha = 1$ one gets the $\Delta$--exponential,
and when $\alpha = 0$ one gets the $\nabla$--exponential, these combined-exponentials can not be really called an exponential function. Indeed, they seem to fail the most important property of an exponential function: they are not a solution of an appropriate initial value problem.


\section{Regular Time Scales}

In the recent paper \cite{MT}, Theorem~\ref{th:1} is used in order to obtain formulas for the $\diam$--derivative of delta and nabla
exponential functions. For that we need assumptions on functions $\rho$ and $\sigma$.
We begin by recalling the notion of \emph{regular time scale} \cite{GGS}.

\begin{definition}\cite{GGS}
A time scale $\T$ is said to be regular if the following two conditions are satisfied simultaneously:

a) $\forall t\in\T, \ \sigma(\rho(t))=t$;

b) $\forall t\in\T, \rho(\sigma(t))=t$.
\end{definition}

\begin{remark}
If $\T$ is a regular time scale, then both operators $\rho$ and $\sigma$ are invertible with
$\sigma^{-1}=\rho$ and $\rho^{-1}=\sigma$.
\end{remark}

The following statement holds:

\begin{proposition}{\rm \cite{GGS}}
\label{prop:GGS}
A time scale $\T$ is regular if and only if the following two conditions are true:

a) The point $\inf\T$ is right-dense and the point $\sup\T$ is left-dense.

b) Each point of $\T^*$ is either two-sided dense or two-sided scattered.
\end{proposition}

Examples of regular time scales include
$\R$, $c \Z$ ($c>0$), $\overline{q^{\Z}}$, $Q_q$, and
$[-\varepsilon,0] \cup \overline{q^{\Z}}$ with $\varepsilon > 0$;
while $\T = [a,b] \cup [c,d]$ is not regular.

\begin{remark}
If $\T$ is a regular time scale, then
$\T^{\kappa}_{\kappa}=\T^{\kappa} = \T_{\kappa} = \T$. Moreover,
$\sigma(\T) = \rho(\T) = \T$.
\end{remark}
We know from Proposition~\ref{prop:GGS} that
for regular time scales $\T$ each $t \in \T^*$
is either two-sided scattered or two-sided dense.
Next proposition gives direct formulas for the $\diam$--derivative of the exponential functions $e_p(\cdot,t_0)$ and
$\hat{e}_p(\cdot,t_0)$.

\begin{proposition}\label{mt:1}{\rm \cite{MT}}
\label{diamexpder}
Let $\T$  be a regular time scale. Assume that
$t,t_0\in\T$ and $p\in \mathcal{R}\cap \mathcal{R}_{\nu}$.
Then,
\begin{gather*}
e_p^{\diam}(t,t_0)=\left[\al p(t)+\frac{(1-\al)p^{\rho}(t)}{1+\nu(t)p^{\rho}(t)}\right]e_p(t,t_0) \, , \\
\hat{e}_p^{\diam}(t,t_0)=\left[(1-\al)p(t)+
\frac{\al p^{\sigma}(t)}{1-\mu(t)p^{\sigma}(t)}\right]\hat{e}_p(t,t_0) \, .
\end{gather*}
\end{proposition}

From Proposition~\ref{mt:1} we can immediately conclude that function $e_p(\cdot,t_0)$ is a solution of the
initial value problem $y^{\diam}(t)=q(t)y(t)$, $y(t_0)=1$,
where $q(t)=\al p(t)+\frac{(1-\al)p^{\rho}(t)}{1+\nu(t)p^{\rho}(t)}$. However, as illustrated by the next example, in general that is not the unique solution. Indeed, for one eigenvalue there may correspond infinitely many eigenvectors.

\begin{example}
\label{ex:in}
Let $\T=c\Z, c>0$, and consider the $\diam$ dynamic equation $y^{\diam}(t)=0$, $y(t_0)=1$, where $\al\in(0,1)$. Then, the constant function $y_1(t)\equiv 1$ is a solution of such initial value problem. Further, it is easily seen that for a given $\al\in(0,1)$ also $y_2(t)=e_q(t,t_0)$, $\ q={-\frac{1}{\al \mu}}$, is a solution:
$\frac{q^{\rho}}{1+q^{\rho}\nu}=
\frac{1}{(1-\al) c}$, and
$$y_2^{\diam}(t)=
\left(\al q+(1-\al)\frac{q^{\rho}}{1+q^{\rho}\nu}\right)e_q(t,t_0)=\left(\frac{-\al}{\al c}+\frac{1-\al}{(1-\al)c}\right)e_q(t,t_0)\equiv 0 \, .$$
We have also that $y_2(t_0)=e_q(t_0,t_0)=1$. Moreover,
each convex combination of $y_1(t)$ and $y_2(t)$ is also a solution of the equation  $y^{\diam}(t)=0$ with condition $y(t_0)=1$. For that, let us see that if  $w(t)=\beta y_1(t)+(1-\beta)y_2(t)=\beta+(1-\beta)e_q(t,t_0)$, then $w^{\diam}(t,t_0)\equiv 0$ and $w(t_0)=\beta+(1-\beta)=1$. Thus, $w(t)$ is also solution to the problem. To have uniqueness of solution, a second boundary or initial condition is needed. Additionally, we see that $e_{-\frac{1}{\al c}}(\rho(t_0),t_0)=\frac{\al}{\al-1}\neq 1$. Hence, with the second condition $y(\rho(t_0))=1$, $y(t)\equiv 1$ is the unique solution.
\end{example}

Example~\ref{ex:in} shows that if we want to have a unique solution for the problem $y^{\diam}(t)=q(t)y(t)$,
$\al\in(0,1)$, then we need to consider two boundary conditions. Only in the cases $\al=1$ and $\al=0$ the second condition disappear. This will be further investigated in Section~\ref{sec:DADE}. Before that
we introduce the idea of decomposition of a regular time scale into a set of time scales that we call \emph{atomic}.

\begin{definition}\label{def:decomp}
Let $\T$ be a regular time scale. Let $I$ be a nonempty finite set of indices. We say that a set $\mathcal{D}=\{\T_i, i \in I\}$, where $\T = \bigcup_{i \in I} \T_i$, forms
an ordered finite partition of the time scale $\T$
if the following holds:
\begin{enumerate}
\item[(i)] for all
$i \in I \backslash \{\max I\}$ one has
$\T_i^* \cap \T_{i+1}^* = \emptyset$ and $\max \T_i=\min \T_{i+1}$ while $\min\T_1=\min\T$, $\max\T_{\max I}=\max \T$;
\item[(ii)] each $\T_i$, $i \in I$,
is a regular time scale;
\item[(iii)] for all $t \in \T_i^*$, $i \in I$, $t$ is either two-sided scattered or two-sided dense.
\item[(iv)] the set $S=\{s_i: s_i=\T_i \cap \T_{i+1}, i=1, \ldots n-1\}$ is such that for any $\eps>0$ we have at least one of the nonempty sets $(s_i-\eps,s_i)\cap\T_i^*$ or $(s_i,s_i+\eps)\cap\T_{i+1}$  consisting of two-sided scattered points only.
\end{enumerate}
An element $\T_i\in \mathcal{D}$ is called an
\emph{atomic time scale}.
If $D=\{\T\}$, then $\T$ is an atomic time scale.
We call the set $S=\{s_i: s_i=\T_i\cap\T_{i+1}, i=1, \ldots n-1\}$ defined in point (iv) the
\emph{set of switching points in the partition $\mathcal{D}$}.
\end{definition}

\begin{remark}
A partition $\mathcal{D}$ of a regular time scale $\T$ can be finite or infinite. In the present study we restrict ourselves to regular time scales $\T$ with a finite partition.
In the case of infinite partitions, we would change the item (i) in Definition~\ref{def:decomp} to:
for all $i \in I \backslash \{\max I\}$ one has
$\T_i^* \cap \T_{i+1}^* = \emptyset$ and $\max \T_i=\min \T_{i+1}$.
\end{remark}

\begin{remark}
Any atomic time scale $\T_i$ in a decomposition of $\T$ does not admit a non-trivial partition: the only partition of
an atomic time scale $\T_i$ is $\mathcal{D} = \{\T_i\}$.
\end{remark}

\begin{remark}
A regular time scale admits a unique partition into a set of atomic time scales because any two such partitions
have the same set of switching points.
\end{remark}

\begin{example}\label{ex:reg}
(a) If $\T = \R_{-} \cup \overline{q^{\Z}}$, then
$\T = \T_1 \cup \T_2$ with $\T_1 := \R_{-} \cup \{0\}$,
$\T_2 := \overline{q^{\Z}}$.
(b) Let $q \in \Z$ with $q > 1$.
If $\T = Q_q=\{-q^k, 0, q^k, k\in\Z\}$,
then $\T = \T_1 \cup \T_2$ with
$\T_1 := \left\{0, -q^k, k\in \Z\right\}$,
$\T_2 :=\left\{0, q^k, k\in \Z\right\}$.
(c) If $\T = \Z$, then $\T = \T_1 := \Z$.
\end{example}

Any regular time scale admits a decomposition into a finite or infinite partition. In our present investigation we only consider regular time scales $\T$ with a finite partition. Example~\ref{ex:inf:part} shows a regular time scale
with an infinite partition.

\begin{example}
\label{ex:inf:part}
Let $\T=\bigcup_{n\in\N} \left(\T_{2n}\cup\T_{2n+1}\right)$, where $\T_{2n}=[2n,2n+1]$ and $\T_{2n+1}=\{2n+1+0.5p, 2n+2-0.5p, p\in \overline{q^{ \Z_{-}\cup\{0\}}} \}$.
Then $\T$ is a regular time scale with an infinite partition.
\end{example}


\section{Diamond-alpha Dynamic Equations}
\label{sec:DADE}

We consider linear equations with $\diam$ derivatives. Our investigations are based on the fact that on a time scale of isolated points, a first order $\diam$ equation becomes a linear equation of second order when written in terms of $\Delta$ and $\nabla$ notions (\textrm{cf.} Examples~\ref{ex:in} and \ref{ex:fib:gn}). Thus,
to have uniqueness, two boundary conditions have to be used. We introduce the idea with a simple example: if $\T=\Z$ and $\al=0.5$, then the Fibonacci sequence is the unique solution to the boundary value problem $y^{\diam}(t)=0.5y(t)$, $y(0)=1$, $y(1)=1$.

\begin{example}
\label{ex:fib:gn}
Let $\T=\N\cup\{0\}$ and $\al=0.5$. We have:
\[2y^{\diam}(t)= y^{\Delta}(t)+y^{\nabla}(t)=y(t+1)-y(t) + y(t)-y(t-1)=y(t+1)-y(t-1).\]
Consider the $\diam$ equation $y^{\diam}(t)=0.5y(t)$.
Then,
$2y^{\diam}(t)=y(t)$, that is, $y(t+1)-y(t-1)=y(t) \ \Leftrightarrow \ y(t+1)=y(t-1)+y(t)$.
To have a unique solution, we need two initial or two boundary conditions. Moreover, if we take $y(0)=1$ and $y(1)=1$, then the Fibonacci sequence is the unique solution to the boundary value problem  $y^{\diamondsuit_{0.5}}(t)=0.5y(t)$, $y(0)=1$, $y(1)=1$.
For $t\in\T$, the solution is the Fibonacci sequence
with the general formula $$y(t)=\frac{1}{\sqrt{5}}\left(\frac{1+\sqrt{5}}{2}\right)^{t+1}-
\frac{1}{\sqrt{5}}\left(\frac{1-\sqrt{5}}{2}\right)^{t+1}\, .$$
\end{example}

\begin{definition}
Let $\T$ be a time scale and  $I\subseteq \T$.
A $\diam$--dynamic equation of the first order on $I$
is an equation of the form $F(t,y(t),y^{\diam}(t))=0$ where $F:I\times \R^2\rightarrow\R$. If function $F$ is linear with respect to $y(\cdot)$ and $y^{\diam}(\cdot)$, then the $\diam$--dynamic equation is said to be linear.
\end{definition}

We consider linear equations of the form  $y^{\diam}(t)=p(t)y(t)+f(t)$, where $f$ and $p$ are given functions on a certain regular time scale $\T$. Under the assumptions of Theorem~\ref{th:1}, we can write that
$y^{\diam}(t)=\alpha y^{\Delta}(t)+(1-\alpha)y^{\nabla}(t)=\alpha y^{\Delta}(t)+(1-\alpha)y^{\Delta\rho}(t)$ or
$y^{\diam}(t)=\alpha y^{\nabla\sigma}(t)+(1-\alpha)y^{\Delta}(t)$. Then, $y^{\diam}(t)=p(t)y(t)+f(t)$ can be rewritten,
respectively, in the following form:
$\alpha y^{\Delta}(t)+(1-\al)y^{\Delta\rho}(t)=p(t)y(t)+f(t)$
or $\alpha y^{\nabla\sigma}(t)+(1-\al)y^{\nabla}(t)=p(t)y(t)+f(t)$.

\begin{definition}
By $L^{\diam}_p$  we denote the operator $L^{\diam}_p: C_{rl}(\T,\R)\rightarrow C_{rl}(\T,\R)$ defined for $t\in\T^{\kappa}_{\kappa}$
by $L^{\diam}_p y(t)=y^{\diam}(t)-p(t)y(t)$.
If $y\in C_{rl}(\T,\R)$ and $L^{\diam}_p y(t)=f(t)$ for all $t\in\T^{\kappa}_{\kappa}$, then we say that $y(\cdot)$ is a solution of $L^{\diam}_p y=f$ on $\T$.
\end{definition}

\begin{proposition}
The operator  $L^{\diam}_p$ is a linear operator on $C_{rl}(\T,\R)$:
$$
L^{\diam}_p(ay_1+by_2)=aL^{\diam}_p y_1 + b L^{\diam}_p y_2
$$
for all $a,b\in\R$ and $y_1,y_2\in C_{rl}(\T,\R)$.
\end{proposition}
\begin{proof}
It follows from the linearity of the $\diam$--derivative.
\end{proof}

\begin{corollary}
If $y_1$ and $y_2$ are solutions of the equation
$L^{\diam}_p y = 0$, then any linear combination of
$y_1$ and $y_2$ is also a solution of $L^{\diam}_p y = 0$, \textrm{i.e.},
$$
L^{\diam}_p y_1 = 0 \text{ and } L^{\diam}_p y_2 = 0 \Rightarrow
L^{\diam}_p\left(a y_1 + b y_2\right) = 0 \quad \forall a, b \in \mathbb{R} \, .
$$
\end{corollary}

\begin{definition}
If $f(t) \equiv 0$, then the equation $L^{\diam}_p y=f$ is reduced to the homogeneous dynamic equation $L^{\diam}_p y = 0$. Otherwise, the equation $L^{\diam}_p y=f(t)$ is called nonhomogeneous.
\end{definition}

The following standard property holds:

\begin{proposition}
The sum of a solution of the homogeneous equation $L^{\diam}_p y = 0$ with a solution of the nonhomogeneous equation $L^{\diam}_p y=f(t)$ gives a solution to the nonhomogeneous equation.
\end{proposition}
\begin{proof}
Let $w,v\in C_{rl}(\T,\R)$. If $L^{\diam}_p w = 0$ and $L^{\diam}_p v=f$, then $L^{\diam}_p (y+w)=f$.
\end{proof}

\begin{theorem}
\label{th:07.07}
Let $\T$ be  an atomic regular time scale,
$t_0\in\T^*$, $y_0$ be a given constant, and $p(\cdot)\in\mathcal{R}\cap\mathcal{R}_{\nu}$.
Then, the two boundary value problems
\begin{equation}
\label{eqdor:1}
\begin{gathered}
L^{\diam}_p y(t)= 0, \quad t\in\T^*, \quad \al\in[0,1],\\
y(t_0)=y_0, \quad \alpha(1-\alpha)y(\rho(t_0))=\alpha(1-\alpha)y_0;
\end{gathered}
\end{equation}
and
\begin{equation}
\label{eqdor:4}
\begin{gathered} \al\mu(s)y^{\Delta\Delta}(s)+\left(1-p^{\sigma}(s)\mu(s)\right)y^{\Delta}(s)
-p^{\sigma}(s)y(s)=0, \quad s\in\T^*,\\
y(\sigma(s_0))=y_0, \quad \alpha(1-\alpha)y(s_0)=\alpha(1-\alpha)y_0;
\end{gathered}
\end{equation}
coincide.
\end{theorem}

\begin{proof}
We do the proof in three parts: first for $\alpha=0$, then for $\alpha=1$, and finally for any $\alpha\in(0,1)$.

Let $\alpha=0$.
Then, problem (\ref{eqdor:1}) has the form $y^{\nabla}(t)=p(t)y(t)$ with one condition $y(t_0)=y_0$ only.
On the other hand, problem (\ref{eqdor:4}) has the form $\left(1-p^{\sigma}(s)\mu(s)\right)y^{\Delta}(s)-p^{\sigma}(s)y(s)=0$
with the initial condition $y(\sigma(s_0))=y_0$.
Denote $\sigma(s)$ by $t$. Then, the dynamic equation of (\ref{eqdor:4}) can be written as $\left(1-p(t)\mu(\rho(t))\right)y^{\Delta}(\rho(t))-p(t)y(\rho(t))=0$. Because $\mu(\rho(t))=\nu(t)$,
$\nu(t)y^{\nabla}(t)=\left\{\begin{array}{lr} 0, & \mbox{if} \ \nu(t)=0\\ y(t)-y(\rho(t)), & \mbox{if}\ \nu(t)\neq 0\end{array}\right.$, and $y^\nabla(t) = y^\Delta\left(\rho(t)\right)$,
we get (\ref{eqdor:4}) in the form $y^{\nabla}(t)=p(t)y(t)$,
$y(t_0)=y_0$.

Consider now $\alpha=1$.
Then, problem (\ref{eqdor:1}) has the form $y^{\Delta}(t)=p(t)y(t)$, $y(t_0)=y_0$.
We transform the dynamic equation of (\ref{eqdor:4}) for cases $\mu(s)=0$ and $\mu(s) \ne 0$ separately.
For $\mu(s)=0$ we get directly the form $y^{\Delta}(s)=p(s)y(s)$, because then $\sigma(s)=s$.
For $\mu(s) \ne 0$ we have that $\mu(s)y^{\Delta\Delta}(s)=y^{\Delta}(\sigma(s))-y^{\Delta}(s)$.
In this case the dynamic equation (\ref{eqdor:4})
can be written as $y^{\Delta}(\sigma(s))-p(\sigma(s))\left(\mu(s)y^{\Delta}(s)-y(s)\right)=0$.
Letting $t=\sigma(s)$, we get $y^{\Delta}(t)=p(t)y(t)$. Moreover, we have also the initial condition $y(t_0)=y(\sigma(s_0))=y_0$.

Let $\al\in(0,1)$, and $s=\rho(t)$. Since $\T$ is a regular time scale, then $t=\sigma(s)$ and
$y^{\diam}(t)=\al y^{\Delta\sigma}(s)+(1-\al)y^{\nabla\sigma}(s)=\al y^{\Delta\sigma}(s)+(1-\al)y^{\Delta}(s)$.
For $\mu(s)=0$ we have also $\nu(s)=0$, and we get the dynamic equation of problem \eqref{eqdor:4} in the form $y^{\Delta}(s)=p(s)y(s)$. In this situation one has $y^{\Delta}(s) = y^{\nabla}(s) = y^{\diam}(s)$ for all $\al\in(0,1)$.
Finally, let us consider points $s\in\T$ which are scattered, \textrm{i.e.}, $\mu(s)\neq 0$ and $\nu(s)\neq 0$, since we have a regular time scale. Then, $\mu(s)y^{\Delta\Delta}(s)=y^{\Delta}(\sigma(s))-y^{\Delta}(s)
=y^{\Delta}(t)-y^{\Delta}((\rho(t))=y^{\Delta}(t)-y^{\nabla}(t)$ and the dynamic equation (\ref{eqdor:4}) simplifies to $\alpha y^{\Delta}(t)+(1-\alpha) y^{\nabla}(t)-p(t)\left(\nu(t)y^{\nabla}(t)-y(\rho(t))\right)=0$.
Thus, $y^{\diam}(t)=p(t)y(t)$.
\end{proof}

\begin{theorem}
\label{thm:mr2}
Let $\T$ be an atomic regular time scale,
$t_0\in\T^*$, $y_0$ be a given constant, and $p(\cdot)\in\mathcal{R}\cap\mathcal{R}_{\nu}$.
Then, the two boundary value problems \eqref{eqdor:1} and \eqref{eqdor:4} have the same  unique solution for each $\al\in[0,1]$.
\end{theorem}

\begin{proof}
For $\alpha=0$ the boundary value problem has the form $y^{\nabla}=p(t)y(t)$, $y(t_0)=y_0$, and
there exists a unique solution for this $\nabla$ initial value problem of first order for all $t\in\T$, because of $\nu$--regressivity of $p$. Similarly,
for $\alpha=1$ the problem has the form $y^{\Delta}=p(t)y(t)$, $y(t_0)=y_0$, and there is a unique solution in $t\in\T$,
because $p \in\mathcal{R}$.
Assume now, up to the end of the proof, that $\al\in(0,1)$.
We consider two cases: $\nu(t) \ne 0$ and $\nu(t) = 0$.
Assume that $\nu(t) \ne 0$. Then, $\mu(s) \ne 0$ and we can write the equation~(\ref{eqdor:4}) in the form
$\al \left[y^{\Delta\sigma}(s)- y^{\Delta}(s)\right]+y^{\Delta}(s)-p^{\sigma}(s)y^{\sigma}(s)=0$.
Moreover, then $\al y^{\Delta}(t)+(1-\al)y^{\Delta}(\rho(t))-p(t)y(t)=0$ and
 $\al y^{\Delta}(t)+(1-\al)y^{\nabla}(t)-p(t)(t)=0$ is in the form of ~(\ref{eqdor:1}).
In this case  equation~(\ref{eqdor:4}) can be written as
\begin{equation*}
y^{\Delta\Delta}(s)+\frac{1-p^{\sigma}(s)\mu(s)}{\al\mu(s)}y^{\Delta}(s)
-\frac{p^{\sigma}(s)}{\al\mu(s)}y(s)=0,
\end{equation*}
which is regressive, thus having a unique solution. Indeed, let  $\tilde{p}(s)=\frac{1-\ps(s)\mu(s)}{\al\mu(s)}$ and $\tilde{q}(s)=-\frac{\ps(s)}{\al\mu(s)}$. Then, for any
$\al\in(0,1)$ one has
$1-\mu(s)\tilde{p}(s)+\mu^2(s)\tilde{q}(s)=1-\frac{1}{\al}
\neq 0$. Finally, consider the case when $\nu(t)=0$. We have
$y^{\Delta}(t)=y^{\nabla}(t)=p(t)y(t)$, $y(t_0) = y_0$, and
we are back to the same situation of $\alpha=0$ or $\alpha=1$.
\end{proof}


\section{Diamond-alpha Exponential Functions}

For a function $y:\T\rightarrow\R$ we define the mapping $Y:\T\rightarrow \R^2$ by setting  $Y(t):=\left[\begin{array}{cc}y(t) & y^{\rho}(t)\end{array}\right]$. Let $t_0,t\in\T$ and $t\geq t_0$. Associated with a boundary value problem
\begin{equation}
\label{eqdor:22.10}
\begin{gathered}
L^{\diam}_p y(t)= 0, \quad t\in\T, \quad \al\in[0,1],\\
y(t_0)=y_0, \quad \alpha(1-\alpha)y(\rho(t_0))=\alpha(1-\alpha)y_0 \, ,
\end{gathered}
\end{equation}
we consider a set of mappings $S(t_0,t):\R^2\rightarrow \R^{2 \times 2}$ for which we can write that $Y(t)=Y(t_0)S(t_0,t)$,  where $y(\cdot)$ is the solution of (\ref{eqdor:22.10}).

Let $\T$ be an atomic regular time scale and let $t\in\T^*$. Then all points in $\T^*$ are two-sided scattered or two-sided dense.
Assume that $\T^*$ consist of scattered points. In that case $\inf \T$ can be right-dense and $\sup\T$ left-dense. For $t\in\T^*$ we can rewrite
$L^{\diam}_p y(t)= 0$ as a second order recurrence equation
using a transition matrix.
If $t_0\in \T$ is right-scattered and $t\geq t_0$, then
\begin{equation}
\label{eq:rec1} Y(t)=Y(t_0)\prod_{k=0}^rA(\sigma^k(t_0))=Y(t_0)A(t_0)A(\sigma(t_0))\ldots A(\sigma^r(t_0)) \, ,
\end{equation}
where $r$ is such that $\sigma^r(t_0)=\rho(t)$ (or $t=\sigma^{r+1}(t_0)$) and matrix $A(t)=\left[\begin{array}{rr}a(t) & 1 \\ b(t) & 0\end{array}\right]$ with
$a(t)=1+\frac{\mu(t)p(t)}{\alpha}-b(t)$,  $b(t)=\frac{(1-\alpha)\mu(t)}{\alpha \nu(t)}$.
The only possibility of $t$ being a left-dense point is that
$t = \sup\T$. In that case $r$ tends to infinity.
For a right-dense $t_0=\inf\T$ we have that
\begin{equation}
\label{eq:rec2} Y(t)=Y(t_0)\left(\prod_{k=1}^{\infty}A^T(\rho^k(t))\right)^T \, .
\end{equation}
Moreover, if $\T$ is a real interval we can also write that
 \begin{equation}\label{eq:rec3} Y(t)=Y(t_0)\exp\left(\int_{t_0}^t p(\tau)d\tau\right)I_{2\times 2} \, .\end{equation}
Hence, we can always write $Y(t)=Y(t_0)S(t_0,t)$, where $S(t_0,t)$ is one of the matrices used in formulas (\ref{eq:rec1})--(\ref{eq:rec3}).
The convergence of the product of matrices depends if $Y(t)$ is defined or not, \textrm{i.e.}, depends on the existence of a solution $y(t)$ to \eqref{eqdor:22.10}.

When we consider a regular time scale $\T$, then we use
its unique partition on atomic time scales.
Let $s_i=\T_i\cap\T_{i+1}$, $t_0\in\T_i, t\in\T_{i+1}$. Then we define $Y(t):=Y(t_0)S(t_0,s_i)S(s_i,t)$, if both $S(t_0,s_i)$ and $S(s_i,t)$ exist.
\begin{remark}
We gave the construction for the forward solution of \eqref{eqdor:22.10}. One can
also consider the backward solution by using the backward formula of recurrence and getting $Y^{\rho}(t)$ from $Y(t)$.
\end{remark}
\begin{definition}
\label{defn:exp:da}
The diamond--alpha exponential function, denoted by $E_{\alpha,p}(\cdot,t_0)$, is defined as the solution
(if it exists) of the diamond--alpha boundary value problem
\begin{equation*}
\begin{gathered}
L^{\diam}_p y(t)= 0, \quad t\in\T, \quad \al\in[0,1],\\
y(t_0)=1, \quad \alpha(1-\alpha)y(\rho(t_0))=\alpha(1-\alpha).
\end{gathered}
\end{equation*}
\end{definition}
\begin{remark}
$E_{\alpha, p}(t,t_0)=\left[\begin{array}{cc} 1 & 1 \end{array}\right]S(t_0,t)\left[\begin{array}{c}1\\0 \end{array}\right]$.
\end{remark}

The following proposition is a direct consequence of the above construction and our Definition~\ref{defn:exp:da}.
\begin{proposition}
Let $\T$ be a regular time scale with a finite partition
and $t_0\in\T$. For all $t\in\T$ we have:

i) if $\al=1$,
$E_{1,p}(t,t_0)=e_p(t,t_0)$;

ii) if  $\al=0$, $E_{0,p}(t,t_0)=\hat{e}_p(t,t_0)$;

iii) for all $\al\in[0,1]$, $E_{\al,0}(t,t_0)\equiv 1$.
\end{proposition}

\begin{example}
Let $\T=\overline{q^{\Z}}$, $q>1$. Then, $\T$ is an atomic time scale. Let $t_0=0$, $p(t)$ be bounded on $\T$,
and $t>0$. Then, $S(t_0,t)=\left(\prod_{k=1}^{\infty}A^T(\rho^k(t))\right)^T$. As $\rho^k(t)$ tends to 0, then
$S(t_0,t)$ exists if $b(t)=\frac{(1-\alpha)q}{\alpha}<1$, \textrm{i.e.}, the diamond--alpha exponential
is well defined if $\frac{q}{q+1}<\alpha<1$.
\end{example}


\begin{example}
Let $\T=c\Z$, $c>0$. Then, $\mu(t)=\nu(t)=c$, $\sigma(t)=t+c$,
and $\rho(t+c)=t$. Let us consider the boundary value problem $y^{\diam}(t)=py(t)$, $y(t_0)=1$,
$\al(\al-1) y(t_0-1)=\al(\al-1)$.
Then, $\tilde{p}(s)=\frac{1-p c}{\al c}$, $\tilde{q}(s)=-\frac{p}{\al c}$, and the dynamic equation
takes the form $y^{\Delta^2}(s)+\tilde{p}(s)y^{\Delta}(s)+\tilde{q}(s)y(s)=0$, which is an equation of constant coefficients and initial conditions $\al(\al-1)y(s_0)=\al(\al-1)$ and $y(\sigma(s_0))=1$.
Consider the problem with $\al=0.5$. Then, the fundamental system of solutions is given by $e_{\lambda_1}(s,s_0)$ and $e_{\lambda_2}(s,s_0)$, where $\lambda_{1,2}=\frac{p c - 1 \pm\sqrt{1+p^2 c^2}}{c}$. To find the exact formula for the solution of the initial value problem we have to calculate $e^{\sigma}_{\lambda}(s,s_0)$ for constant $\lambda$: $e^{\sigma}_{\lambda_{1,2}}(s_0,s_0)
=(1+\mu(s_0)\lambda_{1,2})e_{\lambda_{1,2}}(s_0,s_0)
=1+c\lambda_{1,2}$. The solution has the form
$y(t)=C_1e_{\lambda_1}(s,s_0)+C_2e_{\lambda_2}(s,s_0)$, where constants $C_1$ and $C_2$ satisfy the system $C_1+C_2=1$ and $(1+\lambda_1c)C_1+(1+\lambda_2c)C_2=1$.
As we have that $c(\lambda_1-\lambda_2)=2\sqrt{1+p^2c^2}$, $1+\lambda_1c=p c+\sqrt{1+p^2c^2}$, $1+\lambda_2c=p c-\sqrt{1+p^2c^2}$, then
\[
E_{\frac{1}{2},p}(t,t_0) = \left(\frac{1}{2}
-\frac{p c - 1}{2\sqrt{1+p^2 c^2}}\right)e_{\lambda_1}(t,t_0)
+\left(\frac{1}{2}+\frac{p c - 1}{2\sqrt{1
+ p^2 c^2}}\right)e_{\lambda_2}(t,t_0).
\]
\end{example}


\section*{Acknowledgments}

DM was partially supported by
Bia\l ystok Technical University grant S/WI/1/08;
DT by the R\&D unit CEOC, via FCT and the EC fund FEDER/POCI 2010. This work has flourished during August 2008 while both authors were a Visiting Faculty at the African University of Science and Technology (AUST) in Abuja, Nigeria; and was revised February 2009 also at AUST-Abuja. The good working conditions at AUST Galadimawa Campus are gratefully acknowledged. We would like to thank
an anonymous referee for helpful comments.



\end{document}